\documentclass[twoside]{amsart}
\usepackage{amsmath,amsthm,amssymb}

\theoremstyle{plain}
\newtheorem{theorem}{Theorem}[section]
\newtheorem{proposition}[theorem]{Proposition}
\newtheorem{lemma}[theorem]{Lemma}

\newtheorem{conjecture}[theorem]{Conjecture}
\newtheorem{keylemma}[theorem]{Key-Lemma}

\theoremstyle{remark}
\newtheorem{remark}[theorem]{Remark}
\newtheorem{acknowledgments}{Acknowledgments}

\renewcommand{\k}{{\bf k}}
\newcommand{\dk}{\dot{k}}
\newcommand{\tf}{\tilde f}
\newcommand{\ta}{\tilde a}
\newcommand{\K}{{\bf K}}
\newcommand{\Z}{\mathbb{Z}}
\newcommand{\Q}{\mathbb{Q}}
\newcommand{\QW}{\mathbb{Q}W}
\newcommand{\F}{\mathbb{F}}
\newcommand{\R}{\mathbb{R}}
\newcommand{\N}{\mathbb{N}}
\newcommand{\A}{\mathbb{A}}
\renewcommand{\P}{\mathbb{P}}

\newcommand{\Hd}{ H^1_{\rm dR}   }
\newcommand{\Hc}{ H^1_{\rm crys} }
\newcommand{\HP}{{\rm HP}}
\newcommand{\NP}{{\rm NP}}
\newcommand{\Tr}{{\rm Tr}}
\newcommand{\ord}{ {\rm ord} }

\newcommand{\im}{{\rm Im}}
\newcommand{\res}{{\rm res}}

\newcommand{\cC}{{\mathcal{C}}}
\newcommand{\cU}{{\mathcal{U}}}

\begin{document}

\title{Slope estimates of Artin-Schreier curves}
\author{Jasper Scholten and Hui June Zhu}
\address{
Mathematisch Instituut,
Katholieke Universiteit Nijmegen,
Postbus 9010, 6500 GL Nijmegen.
The Netherlands.
}
\email{scholten@sci.kun.nl}

\address{
Department of mathematics,
University of California,
Berkeley, CA 94720-3840.
The United States.
}
\email{zhu@alum.calberkeley.org}

\date{\today}

\begin{abstract}
Let $f(x)=x^d+a_{d-1}x^{d-1}+\ldots+a_0$ be a polynomial of degree $d$
in $\Q[x]$. For every prime number $p$ coprime to $d$ and $f(x)\in
(\Z_p\cap\Q)[x]$, let $X/\F_p$ be the Artin-Schreier curve defined by
the affine equation $y^p-y=f(x)\bmod p$. Let $\NP_1(X/\F_p)$ be the
first slope of the Newton polygon of $X/\F_p$.  We prove that there is
a Zariski dense subset $\cU$ in the space $\A^d$ of degree-$d$ monic
polynomials over $\Q$ such that for all $f(x)\in\cU$ the following
limit exists and $\lim_{p\rightarrow\infty}\NP_1(X/\F_p)=\frac{1}{d}.$
This is a ``first slope version'' of a conjecture of Wan.

Let $X/\overline\F_p$ be an Artin-Schreier curve defined by the affine
equation $y^p-y=\tf(x)$ where
$\tf(x)=x^d+\ta_{d-1}x^{d-1}+\ldots+\ta_0$.  We prove that if $p>d\geq
2$ then $\NP_1(X/\overline\F_p)\geq
\frac{\lceil\frac{p-1}{d}\rceil}{p-1}$.  If $p>2d\geq 4$, we give a
sufficient condition for the equality to hold.

\end{abstract}

\maketitle

\section{Introduction}\label{section1}

In this paper a curve is a smooth, projective and geometrically
integral algebraic variety of dimension one.
Let $d$ be a natural number.  Let $p$ be a prime coprime to $d$.  Let
$q=p^\nu$ for some natural number $\nu$.  Let $X$ be an {\em Artin-Schreier
curve} over $\F_q$ defined by an affine equation $X: y^p-y
=\tf(x)$, where $\tf(x)\in\F_q[x]$ is of degree $d$.
Then $X$ has genus $g:=\frac{(p-1)(d-1)}{2}$. (See Section 3.)

Write the $L$ function of $X$ over $\F_q$ as
\begin{eqnarray}\label{E:zeta}
\exp\left(\sum_{n=1}^{\infty}(q^n+1-\#X(\F_{q^n}))\frac{T^n}{n}\right)
&=&\frac{1}{P(T)}.
\end{eqnarray}
The denominator $P(T)$ is a polynomial
$1+\sum_{\ell=1}^{2g}b_\ell T^\ell \in 1+T\Z[T]$.
Consider the sequence of points
$(0,0),(1,\frac{\ord_pb_1}{\nu}),(2,\frac{\ord_pb_2}{\nu}),
\ldots,(2g,\frac{\ord_pb_{2g}}{\nu})$ in $\R^2$. (If $b_\ell=0$, define
$\ord_p b_\ell =\infty$.)  The {\em normalized $p$-adic Newton
polygon} of $P(T)$ is defined to be the lower convex hull of this set
of points.  It is called the {\em Newton polygon} of $X/\F_q$,
denoted by $\NP(X/\F_q)$.  Let $\NP_1(X/\F_q)$ denote the first slope
of $\NP(X/\F_q)$, which we call the {\em first slope} of $X/\F_q$.

For any real number $t$ let $\lceil{t}\rceil$ denote the least integer
greater than or equal to $t$ and let $\lfloor{t}\rfloor$ be the biggest
integer less than or equal to $t$.

Let $R$ be a commutative ring with unity.  For any $f(x)\in R[x]$, any
natural numbers $N$ and $r$, we use $[f(x)^N]_r$ to denote the
$x^r$-coefficient of $f(x)^N$.

\begin{theorem}\label{T:box}
Fix $d\geq 2$. Let $X/\F_q$ be an Artin-Schreier curve
of genus $\geq3$ whose
affine equation is given by $y^p-y =\tf(x)$ where
$\tf(x)$ is monic of degree $d$.
\begin{enumerate}
\item[a)] If $p>d$ then
$\NP_1(X/\F_q)\geq \frac{\lceil\frac{p-1}{d}\rceil}{p-1}.$
\item[b)] If $p>2d$ and
$[\tf(x)^{\lceil\frac{p-1}{d}\rceil}]_{p-1}\neq 0$ then
 $\NP_1(X/\F_q)=\frac{\lceil\frac{p-1}{d}\rceil}{p-1}$.
\end{enumerate}
\end{theorem}

Let $\A^d$ be the set of all polynomials
$f(x)= x^d+a_{d-1}x^{d-1}+\ldots + a_0\in\Q[x]$.

\begin{theorem}\label{T:1}
There exists a Zariski dense subset $\cU$ in $\A^d$ such that for all
$f(x)\in \cU$ the Artin-Schreier curve $X/\F_p$ defined by
$y^p-y=f(x)\bmod p$ satisfies
$$\lim_{p\rightarrow\infty}\NP_1(X/\F_p) = \frac{1}{d}.$$
\end{theorem}

If a curve is defined over a perfect field of characteristic $p$
then its Newton polygon is
defined by the ``formal types'' of the $p$-divisible groups associated
to the Jacobian of the curve (see~\cite{Manin} or~\cite{Lenstra-Oort}).
Theorem~\ref{T:box} holds if $\F_q$ is replaced by any perfect field
of characteristic $p$ because its proof remains valid.
It is known that these Newton polygons have integral bending points
and is symmetric in the sense that any line segment of slope $\lambda$
of length $\ell$ occurs in companion with a line segment of slope
$1-\lambda$ of the same length.

Artin-Schreier curves are precisely those degree $p$ abelian covers of
the projective line with the point at infinity totally ramified. So
their $p$-ranks are zero by the Deuring-Shafarevich formula
(see~\cite{Nakajima}).  The $p$-rank is exactly equal to the length of
the slope zero segment of their Newton polygons~(see~\cite{Manin}).
Thus an Artin-Schreier curve has no zero slope.  Suppose $g=1$ or
$2$, then an Artin-Schreier curve $X$ has its first slope equal to
$1/2$.

When $f(x)$ is a monomial then the Frobenius and Verschiebung maps on
the first crystalline cohomology of $X$ have explicit interpretations
(see~\cite{Katz2} and~\cite{Koblitz2}), which enable one to describe
the entire Newton polygon of $X$ explicitly. (Note that classical
literature often refer to this special case as the definition of
Artin-Schreier curve.)

This paper is organized as follows: We recall relevant preliminaries
in Section 2.  Then we develop a method in Section 3 to estimate the
first slopes of Artin-Schreier curves. With some technical preparation
in Sections 4 and 5, Section 6 proves a lower bound for the first
slopes of Artin-Schreier curves, and gives a sufficient condition for
the lower bound to be achieved. We prove Theorem~\ref{T:box} here.
Section 7 detours to prove a similar result by a completely different
method.  Finally in Section 8 we prove Theorem~\ref{T:1}. We will
explain how this is related to the $L$ functions of exponential sums
and a conjecture of Wan.

\begin{acknowledgments}
We thank Daqing Wan for explaining his conjectures and for his
suggestion of Proposition~\ref{P:wan}. We doubly thank him for
numerous stimulating conversations and historical remarks. We
especially thank Bjorn Poonen for his contribution to the proof of
Theorem~\ref{T:2}.
\end{acknowledgments}

\section{Sharp slope estimates}\label{section2}

This section provides fundamental ingredients for our slope estimates
of curves over finite fields.  Note that lemmas we need hold valid
when the base field is perfect of characteristic $p$. However, for simplicity
we constrain ourselves to finite fields in this paper.
Firstly we establish a variation of Katz's
sharp slope estimates in Theorem~\ref{slope}.  Secondly we recall a
method of computing the Verschiebung action on
the first de Rham cohomology of a curve by taking
power series expansions at a rational point. This section essentially
follows~\cite{Katz}, \cite{Katz2}, \cite{Berthelot},
\cite{Berthelot-Ogus} and
\cite{Nygaard:83}. Our approach is particularly inspired by
Nygaard's paper~\cite{Nygaard:83}.

Let $W$ be the ring of Witt vectors over $\F_q$, and $\sigma$ the
absolute Frobenius automorphism of $W$.  Throughout this section we
assume that $X/\F_q$ is a curve of genus $g$ with a rational point.
Let $X/W$ be a smooth and proper lifting of $X$ to $W$, together with
a lifted rational point $P$. The Frobenius endomorphism $F$ (resp.,
Verschiebung endomorphism $V$) are $\sigma$ (resp., $\sigma^{-1}$)
linear maps on the first crystalline cohomology $\Hc(X/W)$ of $X$ with
$FV=VF=p$. It is known that $\Hc(X/W)$ is canonically isomorphic to
the first de Rham cohomology $\Hd(X/W)$ of $X$, and one gets $F$ and $V$
actions on $\Hd(X/W)$.  Thus the pair $(\Hd(X/W),F)$ is a
$\sigma$-$F$-crystal, whereas the pair $(\Hd(X/W),V)$ is a
$\sigma^{-1}$-$V$-crystal.  The Newton polygon of $X/\F_q$ is equal to
the Newton polygon of the crystals $(\Hd(X/W),F)$ and $(\Hd(X/W),V)$
as defined in~\cite{Katz}.

Below we will briefly describe some techniques to approximate slopes
of these crystals.  Let $L$ be the image of $H^0(X,\Omega_{X/W}^1)$ in
$\Hd (X/W)$, and let $M$ be a complement of $L$ such that $\Hd (X/W)=L
\oplus M$. The following lemma is prepared for the proof of
Theorem~\ref{slope}.

\begin{lemma}\label{mazurlemma}
Let notation be as above. Then $L\subset V(L\oplus M)\subset L\oplus
pM$.  If $p^{m-1}$ divides $V^a L$ for some $m>0$ and $a\geq 0$, then
for all $n>a$ we have
\begin{eqnarray*}
V^nL &\subset & p^{m-1}L+p^{m}M.
\end{eqnarray*}
\end{lemma}
\begin{proof}
Recall an equality due to Mazur and Ogus
(see~\cite[Theorem 3]{Mazur2}).
\begin{eqnarray*}
F^{-1}(p(L\oplus M))\otimes\F_q &=&L\otimes\F_q.
\end{eqnarray*}
One easily verifies the following inclusions
$$L\subset F^{-1}(p(L\oplus M))+p(L\oplus M)\subset F^{-1}FV(L\oplus M)
+VF(L\oplus M)\subset V(L\oplus M).$$
The rest follows from~\cite[Lemmas 1.4 and 1.5]{Nygaard:83}.
\end{proof}

\begin{theorem}\label{slope}
Let $\lambda$ be a rational number with $0\leq\lambda\leq\frac{1}{2}$.
Then $\NP_1(X/\F_q)\geq\lambda$ if and only if
\begin{eqnarray*}
p^{\lceil n\lambda\rceil} &\mid& V^{n+g-1}L
\end{eqnarray*}
for all integer $n\geq 1$.
\end{theorem}
\begin{proof}
The main ingredient of the proof is
Katz's sharp slope estimate~\cite[Theorem (1.5)]{Katz}, which
says that $\NP_1(X/\F_p)\geq\lambda$
if and only if $p^{\lceil n\lambda\rceil }\mid V^{n+g}$
for all $n\geq 1$.

Suppose $\NP_1(X/\F_p)\geq\lambda$. Then
$p^{\lceil n\lambda\rceil}\mid V^{n+g}$ for all $n\geq 1$.
By Lemma~\ref{mazurlemma} we have
\begin{eqnarray*}
V^{n+g-1}(L)\quad\subset\quad V^{n+g-1}(\im V)
&=& V^{n+g}(L\oplus M)\quad\subset\quad p^{\lceil n\lambda\rceil}(L\oplus M).
\end{eqnarray*}

Conversely, suppose that
$p^{\lceil n\lambda\rceil} \mid V^{n+g-1} L$
for all $n\geq 1$.
It suffices to show that
$p^{\lceil n\lambda\rceil }\mid V^{n+g}$ for all $n\geq 0$.
For $n=0$ this
statement is trivially true. We proceed by induction on $n$.

Assume that $p^{\lceil (n-1)\lambda\rceil }\mid V^{n+g-1}$.
By Lemma~\ref{mazurlemma}
we have $V(L\oplus M)\subset L\oplus pM$. So
\begin{eqnarray*}
V^{n+g}(L\oplus M)&=&V^{n+g-1}V(L\oplus M)\subset
V^{n+g-1}(L\oplus pM),
\end{eqnarray*}
and
\begin{eqnarray*}
p^{\lceil n\lambda\rceil }\quad \mid\quad p^{\lceil (n-1)\lambda\rceil +1}
&\mid& V^{n+g-1}(pM).
\end{eqnarray*} From the hypothesis
$p^{\lceil n\lambda\rceil }\mid V^{n+g-1} L$,
we have
$p^{\lceil n\lambda\rceil }\mid V^{n+g}$.
\end{proof}

\begin{remark}
Let $n$ and $m$ be any natural integers.
If $p^{m-1}$ divides $V^a L$ for some nonnegative integer $a<n$,
following Lemma~\ref{mazurlemma}, the composition of
$\frac{V^n}{p^{m-1}}$ and reduction to $\F_q$ gives a natural
endomorphism of $L\otimes\F_q$.
This endomorphism of $L\otimes \F_q$ is called a {\em
higher Cartier operator}, denoted by $\cC(m,n)$.  The hypothesis in
the theorem above is equivalent to that $\cC(\lceil n\lambda\rceil,
n+g-1)$ is defined and vanishes for all integer $n\geq 1$.  The
underneath philosophy of our slope estimates is replacing the
traditional Cartier operator by this higher Cartier operator.  We will
not explore this terminology further in this paper.
\end{remark}

Let $\hat X/W$ be the formal completion of $X/W$ at the rational point $P$.
If $x$ is a local parameter of $P$,
then every element of $\Hd(\hat X/W)$ can be represented as
$h(x)\frac{dx}{x}$ for some $h(x)\in xW[[x]]$, and
$F$ and $V$ act as follows:
\begin{eqnarray}
F(h(x)\frac{dx}{x})&=&ph^\sigma(x^p)\frac{dx}{x}\nonumber\\
V(h(x)\frac{dx}{x})&=&h^{\sigma^{-1}}(x^{1/p})\frac{dx}{x}
\quad\mbox{where $x^{m/p}=0$ if $p\nmid m$}\label{EV}.
\end{eqnarray}

Denote the restriction map
$\Hd(X/W)\longrightarrow \Hd(\hat X/W)$ by $\res$.

\begin{lemma}\label{res}
The $F$ and $V$ actions on $\Hd(X/W)$ and $\Hd(\hat X/W)$
commute with the restriction map
\begin{eqnarray*}
\res:\Hd(X/W)&\longrightarrow &\Hd(\hat X/W).
\end{eqnarray*}
Furthermore,
\begin{eqnarray*}
\res^{-1}(p\Hd(\hat X/W)) &=& F(\Hd(X/W)).
\end{eqnarray*}
\end{lemma}
\begin{proof}
The first statement follows from~\cite[Lemma 5.8.2]{Katz2}.
The second is precisely~\cite[Lemma 2.5]{Nygaard:83}.
\end{proof}

This lemma will only be used in the proof of Theorem~\ref{T:induction}.

\section{Slope estimates of Artin-Schreier curves}\label{section3}

Assume that $X$ is an Artin-Schreier curve over $\F_q$ defined by an
affine equation $y^p-y=\tf(x)$ where $\tf(x)=
x^d+\ta_{d-1}x^{d-1}+\ldots + \ta_1x$ and $p\nmid d$.  It is easy to
observe that every Artin-Schreier curve over $\overline\F_p$ can be
written in this form (over some suitable $\F_q$). So $X/\F_q$ has a
rational point at the origin.

Let $p$ be coprime to $d$ and $g\geq 3$.
Take a lifting $X/W$ defined by $y^p-y=f(x)$
where $f(x)=x^d+a_{d-1}x^{d-1}+\ldots+a_1x\in W[x]$ with
$a_\ell\equiv \ta_\ell\bmod p$ for all $\ell$.
So $X/W$ has a rational point at
the origin with a local parameter $x$.  The goal of this
section is to prove Theorem~\ref{T:induction}.  In particular, we
shall prove a highly applicable version in
Key-Lemma~\ref{keylemma}.

For any integer $N>0$ and $0\leq i\leq p-2$ let $C_r(i,N)$ be
the $x^r$-coefficient of the power series expansion
of the function $y^i(py^{p-1}-1)^{p^N-1}$ at the origin $P$:
\begin{eqnarray}\label{E:define-c}
y^i(py^{p-1}-1)^{p^N-1}&=&\sum_{r=0}^{\infty}C_r(i,N)x^r.
\end{eqnarray}

We prepare three lemmas before we start to prove Theorem~\ref{T:induction}.
We shall reserve the range of $i$ and $j$ to be as in
Lemma~\ref{L:1}.

\begin{lemma}\label{L:1}
The curve $X/W$ has genus $(d-1)(p-1)/2$, and for
$p-2\geq i\geq0$, $j\geq1$ and $di+pj\leq(p-1)(d-1)-2+p$
the differential forms
$$\omega_{i,j}:=x^jy^i(py^{p-1}-1)^{-1}\frac{dx}{x}$$
form a basis for $L$.
\end{lemma}
\begin{proof}
For the special fibre $X/\F_q$ this
follows immediately from Proposition VI.4.2 of
\cite{Stichtenoth}.
Let $\QW$ be the field of fractions of $W$.
Consider the generic fibre $X/\QW$.
There are no points $(x,y)$ in $X(\overline\QW)$
with $py^{p-1}-1=f'(x)=0$,
since for such $(x,y)$ one can easily show that $y^p-y$ is not
integral over $W$, and $f(x)$ is integral over $W$.
It follows that the
affine part of $X/\QW$ is nonsingular. The affine ramification points of
the map $X\rightarrow \P^1$ defined by $(x,y)\mapsto x$
are those satisfy $py^{p-1}-1=0$. For each
such $y$ there are exactly $d$ corresponding values of $x$, since
$f'(x)\neq 0$ there. So there are $(p-1)d$ ramification
points on the affine part. The function $y^p-y-f(x)$ in $y$ and its first
two derivatives
have no common zeroes, so all affine ramification
points are of index $2$. Let $e_\infty$ be the ramification index
at $\infty$. By Riemann-Hurwitz, we have
$2g-2=-2p+(p-1)(d-1)+e_\infty-1$. It follows that $g\leq(p-1)(d-1)/2$.
But the genus of the special fibre $X/\F_{q}$ is  $(p-1)(d-1)/2$, hence
the genus of $X/\QW$ is $(p-1)(d-1)/2$, and $e_\infty=p$.

The differential form $\frac{dx}{py^{p-1}-1}=\frac{dy}{f'(x)}$ has no
poles at the affine part.  The form $dx$ only has affine zeroes (of
order 1) at points where the map $(x,y)\mapsto x$ ramifies. At these
points, $py^{p-1}-1=0$, so $\frac{dx}{py^{p-1}-1}$ has no affine
zeroes. The degree of a differential form is $2g-2=(p-1)(d-1)-2$,
hence $\frac{dx}{py^{p-1}-1}$ has a zero of order $(p-1)(d-1)-2$ at
$\infty$. The function $x$ has degree $p$ and no poles at the affine
part. Hence it has a pole of order $p$ at $\infty$. Similarly, $y$ has
a pole of order $d$ at $\infty$. So for $p-2\geq i\geq0$, $j\geq1$ and
$di+pj\leq(p-1)(d-1)-2+p$ the form $\omega_{i,j}$ is in $L$. From the
assumption that $d$ and $p$ are coprime it follows that for $i$ and $j$ in
this range the $\omega_{i,j}$ have zeroes of different order at
$\infty$, hence they are independent.  From \cite{Stichtenoth},
Proposition VI.4.2 (h), it follows that the reduction of these
differential forms modulo $p$ form a basis for
$H^0(\Omega^1_{X/\F_{q}})$, hence the $\omega_{i,j}$ form a basis for
$L$.
\end{proof}

\begin{lemma}\label{L:warm-up}
Let $m$ be a positive integer. If $p\nmid m$ then
$x^m(py^{p-1}-1)^{-1}\frac{dx}{x}\equiv 0\bmod p$ in $\Hd(\hat X/W)$.
\end{lemma}
\begin{proof}
If $p\nmid m$ then
$x^m(py^{p-1}-1)^{-1}\frac{dx}{x}\equiv-d\left(\frac{x^m}{m}\right)\bmod
p$, which is zero in $\Hd(\hat X/W)$.
\end{proof}

\begin{lemma}\label{indN}
For all nonnegative integers $a$ and $r$ we have
\begin{eqnarray*}
C_r(i,N+a)&\equiv &C_r(i,N)\bmod p^{N+1}.
\end{eqnarray*}
\end{lemma}
\begin{proof}
It is easy to see that
$\binom{p^N}{\ell}\equiv 0\bmod p^{N+1-\ell}$
if $N+1\geq \ell\geq 1$.
Thus
\begin{eqnarray*}
(1-py^{p-1})^{p^N}
\quad =\quad
\sum_{\ell=0}^{p^N}\binom{p^N}{l}(-py^{p-1})^\ell
\quad \equiv \quad  1\bmod p^{N+1}.
\end{eqnarray*}
Therefore, we have
\begin{eqnarray*}
y^i(py^{p-1}-1)^{p^{N+a}-1}
&=&
y^i(py^{p-1}-1)^{p^N-1}
(1-py^{p-1})^{p^N(p^a-1)}\\
&\equiv &
y^i(py^{p-1}-1)^{p^N-1} \bmod p^{N+1}.
\end{eqnarray*}
This proves the lemma.
\end{proof}

\begin{theorem}\label{T:induction}
Let $\lambda$ be a rational number with $0\leq \lambda\leq\frac{1}{2}$.
Suppose there exists an integer $n_0$ such that
\begin{enumerate}
\item[(i)] for all $m\geq 1$ and $1\leq n<n_0$ we have
\begin{eqnarray*}
\ord_p(C_{mp^{n+g-1}-j}(i,n+g-2))\geq \lceil n\lambda\rceil;
\end{eqnarray*}
\item[(ii)] for all $m\geq 2$ we have
\begin{eqnarray*}
\ord_p(C_{mp^{n_0+g-1}-j}(i,n_0+g-2))\geq \lceil n_0\lambda\rceil.
\end{eqnarray*}
\end{enumerate}
Then
$$
\left\{
\begin{array}{ll}
p^{\lceil{n\lambda}\rceil}   \mid V^{n+g-1}(\omega_{i,j})&
\mbox{if $n< n_0$}; \\
p^{\lceil n_0\lambda\rceil-1}\mid V^{n_0+g-1}(\omega_{i,j})&
\mbox{if $n=n_0$}.
\end{array}
\right.
$$
Furthermore, we have
\begin{eqnarray*}
V^{n_0+g-1}(\omega_{i,j})&=&
C_{p^{n_0+g-1}-j}^{\sigma^{-(n_0+g-1)}}(i, n_0+g-2)(\omega_{0,1}).
\end{eqnarray*}
\end{theorem}
\begin{proof}
We will prove the first part by induction.
Suppose $1\leq n\leq n_0$ and
\begin{eqnarray}\label{eq:indhyp}
p^{\lceil{(n-1)\lambda}\rceil}&\mid& V^{n+g-2}(\omega_{i,j}).
\end{eqnarray}
Note that this this trivially true if $n=1$.

Write $h(x):=\res(py^{p-1}-1)^{-1}$.
By~\cite[Lemma 2.2]{Nygaard:83}, we have
$$
h(x)^{p^{n+g-2}}
=h^{\sigma^{n+g-2}}(x^{p^{n+g-2}})+
ph_1^{\sigma^{n+g-3}}
(x^{p^{n+g-3}})+\ldots+p^{n+g-2}h_{n+g-2}(x).$$
for some power series
$h_1(x),h_2(x),\ldots,h_{g+n-2}(x) \in W[[x]]$.
Thus the power series expansion of $\omega_{ij}$ is
\begin{eqnarray*}
\res(\omega_{ij})
&=&\res\left(x^jy^i(py^{p-1}-1)^{-1}\frac{dx}{x}\right)\\
&=&\res\left(x^jy^i(py^{p-1}-1)^{p^{n+g-2}-1}h(x)^{p^{n+g-2}}
\frac{dx}{x}\right)\\
&=&\sum_{r=0}^\infty C_r(i,n+g-2)x^{r+j}h^{\sigma^{n+g-2}}(x^{p^{n+g-2}})
\frac{dx}{x}\\
&&+p\sum_{r=0}^\infty C_r(i,n+g-2)x^{r+j}h_1^{\sigma^{n+g-3}}
  (x^{p^{n+g-3}})\frac{dx}{x}\\
&&+\cdots\\
&&+p^{n+g-2}\sum_{r=0}^\infty C_r(i,n+g-2)x^{r+j}
h_{n+g-2}(x)\frac{dx}{x}.
\end{eqnarray*}

Apply the $V^{g+n-2}$ action on the above first differential form. Since
the $V$-action commutes with the restriction map (by Lemma~\ref{res}),
we have
\begin{eqnarray}\label{Vaction}
&&\res(V^{n+g-2}\omega_{ij})
=\sum_{m=1}^\infty C^{\sigma^{-(n+g-2)}}_{mp^{n+g-2}-j}(i,n+g-2)
x^mh(x)\frac{dx}{x}\\\nonumber
&&\phantom{\res(V^{n+g-2}\omega_{ij})}
+p\sum_{m=1}^\infty
   C^{\sigma^{-(n+g-3)}}_{mp^{n+g-3}-j}(i,n+g-2)
   V\left(x^mh_1(x)\frac{dx}{x}\right)\\\nonumber
&&\phantom{\res(V^{n+g-2}\omega_{ij})}
+p^2\sum_{m=1}^\infty
   C^{\sigma^{-(n+g-4)}}_{mp^{n+g-4}-j}(i,n+g-2)
   V^2\left(x^mh_2(x)\frac{dx}{x}\right)\\\nonumber
&&\phantom{\res(V^{n+g-2}\omega_{ij})}
+\cdots\\\nonumber
&&\phantom{\res(V^{n+g-2}\omega_{ij})}
+p^{\lceil n\lambda\rceil-1}\sum_{m=1}^{\infty}
C^{\sigma^{-(n+g-1-\lceil n\lambda\rceil)}}
_{i,mp^{n+g-1-\lceil n\lambda\rceil}-j}(n+g-2)
\\\nonumber
&&\phantom{\res(V^{n+g-2}\omega_{ij})}
\phantom{+p^{\lceil n\lambda\rceil-1}}
\cdot V^{\lceil n\lambda\rceil-1}
 \left(x^mh_{\lceil n\lambda\rceil -1}(x)\frac{dx}{x}\right)\\\nonumber
&&\phantom{\res(V^{n+g-2}\omega_{ij})}
+p^{\lceil n\lambda\rceil}\beta,
\end{eqnarray}
for some $\beta\in\Hd(\hat X/W)$.

By the hypothesis,
$p^{\lceil n\lambda\rceil-1}$ divides $C_{mp^{n+g-2}-j}(i,n+g-3)$.
For all $m\geq1$, by Lemma~\ref{indN},
\begin{eqnarray}\label{eq:congr1}
p^{\lceil n\lambda\rceil-1}&\mid&C_{mp^{n+g-2}-j}(i,n+g-2).
\end{eqnarray}
For $m$ coprime to $p$ it follows from Lemma~\ref{L:warm-up} that
$p$ divides $x^mh(x)\frac{dx}{x}$. Thus
\begin{eqnarray*}
p^{\lceil n\lambda\rceil}&\mid& C_{mp^{n+g-2}-j}(i,n+g-2)x^mh(x)\frac{dx}{x}.
\end{eqnarray*}
Otherwise, except possibly when $n=n_0$ and $m=p$, we have
\begin{eqnarray*}
p^{\lceil n\lambda\rceil}&\mid& C_{(\frac{m}{p})p^{n+g-1}-j}(i,n+g-2).
\end{eqnarray*}
Therefore,
\begin{eqnarray}\label{eq:congr2}
&&\sum_{m=1}^\infty C^{\sigma^{-(n+g-2)}}_{mp^{n+g-2}-j}(i,n+g-2)
x^mh(x)\frac{dx}{x}\\\nonumber
&\equiv&
\sum_{m'=1}^\infty C^{\sigma^{-(n+g-2)}}_{m'p^{n+g-1}-j}(i,n+g-2)
x^{pm'}h(x)\frac{dx}{x}\\\nonumber
&\equiv&\left\{
\begin{array}{ll}
0\bmod p^{\lceil n\lambda\rceil}&{\rm if\ }n<n_0\\
C^{\sigma^{-(n_0+g-2)}}_{p^{n_0+g-1}-j}(i,n_0+g-2)
  x^ph(x)\frac{dx}{x}
\bmod p^{\lceil n\lambda\rceil}&
{\rm if\ }n=n_0
\end{array}
\right.
\end{eqnarray}
For all integer $\ell\geq 1$, by the hypothesis of the theorem,
we obtain
\begin{eqnarray*}
\ord_p(C_{mp^{n+g-\ell-2}-j}(i,n+g-\ell-3))
\geq \lceil (n-\ell-1)\lambda \rceil
\geq \lceil n\lambda \rceil - \ell.
\end{eqnarray*}
So, by Lemma~\ref{indN}, we have
$\ord_p(C_{mp^{n+g-\ell-2}-j}(i,n+g-2))\geq\lceil n\lambda\rceil-\ell$.
So $p^{\lceil n\lambda\rceil}$ divides every sum of $(\ref{Vaction})$
except possibly the one on the first line.
Combining this information with $(\ref{eq:indhyp})$, $(\ref{eq:congr1})$
and $(\ref{eq:congr2})$ yields for all $n<n_0$
\begin{eqnarray*}
\res\left(\frac{V^{n+g-2}(\omega_{i,j})}
{p^{\lceil n\lambda\rceil-1}}\right)&\in& p\Hd(\hat X/W)
\end{eqnarray*}
Hence for such $n$ Lemma~\ref{res} implies
\begin{eqnarray*}
\frac{V^{n+g-2}(\omega_{i,j})}
{p^{\lceil n\lambda\rceil-1}}&\in& F(\Hd(X/W))
\end{eqnarray*}
so
\begin{eqnarray*}
\frac{V^{n+g-1}(\omega_{i,j})}
{p^{\lceil n\lambda\rceil-1}} &\in&
VF(\Hd(X/W))=p\Hd(X/W),
\end{eqnarray*}
which proves the induction hypothesis.
If $n=n_0$ then the above implies
\begin{eqnarray*}
\res\left(\frac{V^{n_0+g-2}(\omega_{i,j})}
{p^{\lceil n_0\lambda\rceil-1}}\right)
-\frac{1}{p^{\lceil n_0\lambda\rceil-1}}
C^{\sigma^{-(n_0+g-2)}}_{p^{n_0+g-1}-j}(i,n_0+g-2)
x^{p}h(x)\frac{dx}{x}
\end{eqnarray*}
lies in $p\Hd(\hat X/W).$
Lemma~\ref{res} implies
\begin{eqnarray*}
\frac{V^{n_0+g-2}(\omega_{i,j})}
{p^{\lceil n_0\lambda\rceil-1}}
-  \frac{1}{p^{\lceil n_0\lambda\rceil-1}}
C^{\sigma^{-(n_0+g-2)}}_{p^{n_0+g-1}-j}(i,n_0+g-2)
\omega_{0,p}
\end{eqnarray*}
lies in $F(\Hd(X/W))$. Hence,
\begin{eqnarray*}
 \frac{V^{n_0+g-1}(\omega_{i,j})}
{p^{\lceil n_0\lambda\rceil-1}}
-  \frac{1}{p^{\lceil n_0\lambda\rceil-1}}
C^{\sigma^{-(n_0+g-1)}}_{p^{n_0+g-1}-j}(i,n_0+g-2)
V(\omega_{0,p})
\end{eqnarray*}
lies in $VF\Hd(X/W)=p\Hd(X/W)$.
Now the theorem follows from $V(\omega_{0,p})\equiv\omega_{0,1}\bmod p$.
\end{proof}

We summarize everything we need in the key lemma below.

\begin{keylemma}\label{keylemma}
Let $\lambda$ be a rational number with $0\leq\lambda\leq\frac{1}{2}$.
\begin{enumerate}
\item[i)] If for all $i,j$ within the range of Lemma~\ref{L:1},
and for all $m\geq1$, $n\geq1$ we have
 \begin{eqnarray*}
 \ord_p(C_{mp^{n+g-1}-j}(i,n+g-2))&\geq&\lceil n\lambda\rceil
 \end{eqnarray*}
 then $$\NP_1(X/\F_q)\geq\lambda.$$
\item[ii)] Let $i,j$ be within the range.
	\begin{enumerate}
\item[a)] Let $n_0\geq1$. Suppose that for all
 $m\geq1$ and $1\leq n<n_0$ we have
 \begin{eqnarray*}
 \ord_p(C_{mp^{n+g-1}-j}(i,n+g-2))&\geq&\lceil n\lambda\rceil;
 \end{eqnarray*}
\item[b)] suppose that for all $m\geq2$ we have
 \begin{eqnarray*}
 \ord_p(C_{mp^{n_0+g-1}-j}(i,n_0+g-2))&\geq&\lceil n_0\lambda\rceil;
 \end{eqnarray*}
\item[c)] suppose
 \begin{eqnarray*}
 \ord_p(C_{p^{n_0+g-1}-j}(i,n_0+g-2))&<& \lceil n_0\lambda\rceil;
 \end{eqnarray*}
	\end{enumerate}
Then
 \begin{eqnarray*}
 \NP_1(X/\F_q)&<&\lambda.
 \end{eqnarray*}
\end{enumerate}
\end{keylemma}
\begin{proof}
i) The hypotheses in Theorem \ref{T:induction} are satisfied for all
positive integers $n_0$ and for all possible $i$ and $j$.
Thus our statement follows from Theorem~\ref{slope}.

ii)
If $\NP_1(X/\F_q)\geq \lambda$ then $p^{\lceil n_0\lambda\rceil} \mid
V^{n_0+g-1}(\omega_{i,j})$ for all $i,j$ in the range of Lemma~\ref{L:1}
by Theorem~\ref{slope}.
This implies that for the particular $i,j$
satisfying the hypothesis of Theorem \ref{T:induction} we have
$$\ord_p(C_{p^{n_0+g-1}-j}(i,n_0+g-2))
\geq \lceil n_0\lambda\rceil.$$
This proves the Key-Lemma.
\end{proof}

\section{$p$-adic behavior of coefficients of power series}
\label{section4}

In this section we study the $p$-adic behavior of
coefficients of two power series.

To make this paper as self-contained as possible, we recall
{\em Lagrange inversion formula} from mathematical analysis
~\cite[IX, \S~189]{Goursat}.
Let $z$ and $y$ be two functions such that
$y=z\mu(y)$ for some function $\mu(y)$
which can be developed into a power series in $y$. Then
the power series expansion of any function $h(y)$ in $z$ is
\begin{eqnarray}\label{E:Lagrange}
h(y)&=&
\sum_{k_1=1}^{\infty}
\frac{1}{k_1!}\left(\left.\left(\mu(y)^{k_1} h'(y)\right)^{(k_1-1)}
\right|_{y=0}
\right) z^{k_1},
\end{eqnarray}
where the upper corner $^{(k_1-1)}$ denotes the $(k_1-1)$-th derivative
and $h'(y)$ denotes the first derivative of $h(y)$ in terms of $y$.

\begin{lemma}\label{L:4.1}
Let $a>0$ and let $y\in W[[z]]$ be a power series that satisfies
$y^p-y=z$ and $y(0)=0$. Then
$$y^a=\sum_{k_1=1}^\infty D_{k_1}(a)z^{k_1}$$
where $D_{k_1}(a)=0$ if $k_1\not\equiv a\bmod p-1$; otherwise,
$$D_{k_1}(a)=(-1)^{a+\frac{k_1-a}{p-1}}
\frac{a\,\left(k_1+\frac{k_1-a}{p-1}-1\right)!}
{k_1!\left(\frac{k_1-a}{p-1}\right)!}.$$
\end{lemma}
\begin{proof}
Note that $y=z(y^{p-1}-1)^{-1}$.
Apply (\ref{E:Lagrange}) to this equation, we get
\begin{eqnarray}\label{E:y_t}
y^a = \sum_{k_1=1}^{\infty}
\frac{a}{k_1!}\left(\left.\left((y^{p-1}-1)^{-k_1}y^{a-1}\right)^{(k_1-1)}
\right|_{y=0}
\right)z^{k_1}.
\end{eqnarray}
We have
$$\left.\left(y^{p-1}-1)^{-k_1}y^{a-1}\right)^{(k_1-1)}\right|_{y=0}=
\left.\left(\sum_{\ell=0}^\infty
(-1)^{(p-1)\ell+k_1}\binom{-k_1}
{\ell}y^{(p-1)\ell+a-1}\right)^{(k_1-1)}\right|_{y=0}.$$
Clearly, this is 0 if $k_1\not\equiv a\bmod p-1$; otherwise, it is
equal to
$$(-1)^a(k_1-1)!\binom{-k_1}{\frac{k_1-a}{p-1}}.$$
Plugging this into (\ref{E:y_t}) yields the desired
value for $D_{k_1}(a)$.
\end{proof}

For any natural numbers $k_1$ and $a$,
we will keep the notation $D_{k_1}(a)$ as defined in Lemma~\ref{L:4.1}.
We also define
$D_{k_1}(a)=1$ if $a=k_1=0$ and $D_{k_1}(a)=0$ if only one of $k_1$ and
$a$ is 0.
For any integer $k\geq0$ denote by $s_p(k)$ the sum of all digits in
the ``base $p$'' expansion of $k$.

\begin{lemma}\label{L:4.2}
If $a>0$ and $k_1\equiv a\bmod p-1$,
write $a=i+\ell(p-1)$ with integers $\ell$ and $1\leq i\leq p-1$,
then
$$\left\{\begin{array}{ll}
\ord_p(D_{k_1}(a))=\displaystyle
      \frac{s_p({k_1})-i}{p-1}&\mbox{if $\ell=0$};\\
\ord_p(D_{k_1}(a))\geq\displaystyle
      \frac{s_p({k_1})-i}{p-1}-(\ell-1)&\mbox{if $\ell\geq 1$}.
\end{array}\right.$$
\end{lemma}

\begin{proof}
Let
$k_1\equiv a\bmod p-1$.  Using the well-known identity
$(p-1)\ord_p(k!)=k-s_p(k)$ for all natural number $k$, one gets that
\begin{eqnarray}\label{eq:1}
&&\quad\ord_p(D_{k_1}(a))\\
&&\phantom{\ord_p}=\ord_p(a)+\frac{1}{p-1}\left(
s_p({k_1})+s_p\left(\frac{{k_1}-a}{p-1}\right)-1
-s_p\left(a-1+\frac{{k_1}-a}{p-1}p\right)\right).\nonumber
\end{eqnarray}
If $\ell=0$ then
$$s_p\left(a-1+\frac{k_1-a}{p-1}p\right)
= i-1+s_p\left(\frac{k_1-a}{p-1}\right).$$
If $\ell=1$ then
$$s_p\left(a-1+\frac{k_1-a}{p-1}p\right)
\leq
(p-1)\ord_p(a)+i-1+s_p\left(\frac{{k_1}-a}{p-1}\right).$$
If $\ell>1$ then
\begin{eqnarray*}
s_p\left(a-1+\frac{k_1-a}{p-1}p\right)
&\leq &i-1+s_p(\ell(p-1))+s_p\left(\frac{{k_1}-a}{p-1}\right)\\
&\leq & i-1+(\ell-1)(p-1)+s_p\left(\frac{{k_1}-a}{p-1}\right).
\end{eqnarray*}
Substitute these back in $(\ref{eq:1})$, we obtain
the desired (in)equalities.
\end{proof}

Fix two integers $N>0$ and $0\leq i\leq p-2$.
Let $y\in W[[z]]$ still be the power series satisfying
$y^p-y=z$. Define coefficients $E_{k_1}(i,N)$ by
\begin{eqnarray*}
y^i(py^{p-1}-1)^{p^N-1} &=&\sum_{k_1=0}^\infty E_{k_1}(i,N)z^{k_1}.
\end{eqnarray*}
For any integer $r\geq0$ let $\K_r$ denote the set of
transposes $\k=\ ^t\!(k_1,\ldots,k_d)$ of $d$-tuple integers
with $k_1\geq k_2\geq\ldots\geq k_d\geq 0$ and
$\sum_{\ell=1}^{d} k_\ell=r$.
We define
\begin{eqnarray*}
s_p(\k)&:=&s_p(k_1-k_2)+\ldots+s_p(k_{d-1}-k_d)+s_p(k_d).
\end{eqnarray*}
Note that from the definition of the coefficients $C_r(i,N)$ in
(\ref{E:define-c}) we find
$$\sum_{r=0}^\infty C_r(i,N)x^r=\sum_{k_1=0}^\infty
E_{k_1}(i,N)f(x)^{k_1}.$$
Expanding the powers of $f(x)$ yields
\begin{eqnarray}\label{E:C-E}
C_r(i,N)=\sum_{\k\in\K_r}E_{k_1}(i,N)\prod_{\ell=1}^{d-1}
\binom{k_\ell}{k_{\ell+1}}a_{\ell}^{k_\ell-k_{\ell+1}}.
\end{eqnarray}

\begin{lemma}\label{L:4.3}
Let $\k=\ ^t\!(k_1,\ldots,k_d)\in\K_r$.
If ${k_1}\not\equiv i\bmod p-1$ then
$E_{k_1}(i,N)=0$. If ${k_1}\equiv i\bmod p-1$ then
\begin{eqnarray*}
\ord_p(E_{k_1}(i,N))&=&
\frac{s_p(k_1)-i}{p-1},\\
\ord_p\left(E_{k_1}(i,N)\prod_{\ell=1}^{d-1}\binom{k_\ell}
{k_{\ell+1}}\right)&=&
\frac{s_p(\k)-i}{p-1}.
\end{eqnarray*}
\end{lemma}
\begin{proof}
Take the identity
\begin{eqnarray}\label{E:sum}
y^i(py^{p-1}-1)^{p^N-1} &=& \sum_{\ell=0}^{p^N-1}(-1)^{p^N-1-\ell}
\binom{p^N-1}{\ell}p^\ell y^{i+\ell (p-1)}.
\end{eqnarray}
Substitute the power series expansion of
$y^{i+\ell(p-1)}$ in
(\ref{E:sum}); we get
\begin{eqnarray}\label{E:ck}
E_{k_1}(i,N)&=&
\sum_{\ell=0}^{p^N-1}(-1)^{p^N-1-\ell}\binom{p^N-1}{\ell}
D_{k_1}(i+\ell(p-1))p^\ell.
\end{eqnarray}

If ${k_1}\not\equiv i \bmod p-1$ then $D_{k_1}(i+\ell(p-1))=0$
by Lemma~\ref{L:4.1};
hence $E_{k_1}(i,N)=0$. This prove the first part of the lemma.

If $k_1=i=0$ then $E_{k_1}(i,N)=(-1)^{p-1}$ and $s_p(k_1)-i=0$.
If $i=0$, $k_1>0$ and $i\equiv k_1\bmod p-1$  then, by Lemma~\ref{L:4.2},
the term with minimal valuation in (\ref{E:ck}) occurs
at $\ell=1$. We have
$$\ord_p(E_{k_1}(i,N))=1+\ord_p(D_{k_1}(p-1))=
1+\frac{s_p({k_1})-(p-1)}{p-1}.$$
If $i>0$ and $k_1\equiv i\bmod p-1$ then the term with minimal
valuation in (\ref{E:ck}) occurs at
$\ell=0$. We have
$$\ord_p(E_{k_1}(i,N))=\ord_p(D_{k_1}(i))=
\frac{s_p({k_1})-i}{p-1}.$$
This implies the second assertion.

By $\ord_p(k!)=\frac{k-s_p(k)}{p-1}$ we have
that
$\ord_{p}\binom{k_\ell}{k_{\ell+1}}=
\frac{s_p(k_\ell)+s_p(k_{\ell}-k_{\ell+1})-s_p(k_\ell)}{p-1}$.
Thus
$$\ord_p(\prod_{\ell=1}^{d-1}\binom{k_{\ell}}{k_{\ell+1}})
=\frac{s_p(\k)-s_p(k_1)}{p-1}.$$
So the third assertion follows from this equality and the second
assertion.
\end{proof}

\section{$p$-adic behavior of $C_r(i,N)$}\label{section5}

To apply Theorem \ref{T:induction} one needs to have in hands
an efficient formula of the $p$-adic valuations of the coefficients
in~(\ref{E:define-c}). This formula is in
Lemma~\ref{L:pbox3}, which is prepared for Section~\ref{section6}.

Let $\k=\ ^t\!(k_1,\ldots,k_d)\in\K_r$.
For $1\leq \ell\leq d$, let
$k_\ell=\sum_{v\geq 0}k_{\ell,v}p^v$ be the ``base $p$'' expansion of
$k_\ell$, we introduce a {\em dot representation}
$$\dk_\ell:= [\ldots, \dk_{\ell,2},\dk_{\ell,1}, \dk_{\ell,0}]$$
in the following way: for $\ell=d$, let $\dk_{d,v}=k_{d,v}$
for all $v\geq 0$; for $1\leq \ell<d$, it is defined inductively by
$$\dk_{\ell-1,v}:=\dk_{\ell,v}+ \mbox{$p^v$-coefficient in the ``base $p$''
expansion of $(k_{\ell-1}-k_{\ell})$},$$
for all $v\geq 0$.
It can be verified that
$k_\ell=\sum_{v\geq 0}\dk_{\ell,v}p^v$ for $1\leq \ell\leq d$.  Since
$k_\ell\geq k_{\ell+1}$ we have $\dk_{\ell-1,v}\geq \dk_{\ell,v}$ for
all $v$. It is not hard to observe
$$s_p(\k)=\sum_{v\geq 0}\left(
\sum_{\ell=1}^{d-1}(\dk_{\ell,v}-\dk_{\ell+1,v})+\dk_{d,v}\right)
= \sum_{v\geq 0}\dk_{1,v}.$$
For any natural number $a$, define a subset of $\K_r$ as follows
\begin{eqnarray*}
\K_r^a&:=&\{\k\in\K_r\ |\ \dk_{\ell,v}=0 {\rm\ for
\ }v\geq a,\ 1\leq\ell\leq d\}.
\end{eqnarray*}
More explicitly $\K_r^a$ consists of all $\k\in\K_r$ with
$\dk_\ell =[\ldots,0,\dk_{\ell,a-1},\ldots,
\dk_{\ell,1},\dk_{\ell,0}]$ for all $1\leq \ell\leq d$.
Then we have an obvious filtration
$\K_r^1\subseteq\ldots
\subset \K_r^{a-1}\subseteq\K_r^{a}\subseteq\ldots\subseteq \K_r$.

\begin{lemma}\label{L:pbox1}
Let $p>d$. Let $1\leq j\leq p-1$, let $a, m,n\geq1$ and
$r=mp^a-j$. If $\k\in\K_r^a$, then
\begin{eqnarray*}
s_p(\k)&\geq& \left\lfloor\frac{(m-1)p}{d}\right\rfloor+
(a-1)\left\lceil\frac{p-1}{d}\right\rceil+\left\lceil\frac{p-j}{d}\right\rceil.
\end{eqnarray*}
If $p>2d$, $m=1$ and the equality holds, then
$$\dk_1 =
[\ldots,
0,
\overbrace{
\lceil\frac{p-1}{d}\rceil,
\ldots,
\lceil\frac{p-1}{d}\rceil,
\lceil\frac{p-j}{d}\rceil}^{a}
].$$
\end{lemma}
\begin{proof}
We will prove this lemma by induction on $a$. Let $\k\in\K_r^a$.
Suppose $a=1$. Note that
all real numbers $c_1$ and $c_2$ satisfy
$\lceil c_1+c_2\rceil\geq\lfloor c_1\rfloor+\lceil c_2\rceil$.
So we have
$$s_p(\k)=\dk_{1,0}\geq\left\lceil\frac{r}{d}\right\rceil\geq
\left\lfloor\frac{(m-1)p}{d}\right\rfloor
+\left\lceil\frac{p-j}{d}\right\rceil.$$
Suppose $m=1$ and the equality holds.
It reads $\dk_{1,0}=\lceil\frac{p-j}{d}\rceil$.

Now suppose $a\geq 2$.
Let $\k\in\K_r^a$.
Let $k'_\ell:=\sum_{v=0}^{a-2}\dk_{\ell,v}p^v$ for all $1\leq\ell\leq d$.
One can find a natural number $m'$ such that
$m'p^{a-1}-j=\sum_{\ell=1}^d k_\ell'.$
Then $\k':={}^t\!(k'_1,\ldots,k'_d)\in \K_{m'p^{a-1}-j}^{a-1}$.
We have $(m-1)p+p-1=(m'-1)+\sum_{\ell=1}^d\dk_{\ell,a-1}\leq (m'-1)
+d\dk_{1,a-1}$.
So
\begin{eqnarray}\label{eq:pbox1}
\dk_{1,a-1}&\geq &\left\lceil\frac{(m-1)p-(m'-1)+p-1}{d}\right\rceil\\
\nonumber
&\geq &\left\lfloor\frac{(m-1)p}{d}\right\rfloor-\left\lceil
\frac{m'-1}{d}\right\rceil
+\left\lceil\frac{p-1}{d}\right\rceil.
\end{eqnarray}
On the other hand, by induction hypothesis on $\k'\in\K_{m'p^{a-1}-j}^{a-1}$,
one has
\begin{eqnarray}\label{eq:pbox2}
\sum_{v=0}^{a-2} \dk_{1,v}\geq
\left\lfloor\frac{(m'-1)p}{d}\right\rfloor
+(a-2)\left\lceil\frac{p-1}{d}\right\rceil
+\left\lceil\frac{p-j}{d}\right\rceil.
\end{eqnarray}
Combining (\ref{eq:pbox1}) and (\ref{eq:pbox2}), one gets
\begin{eqnarray}\label{eq:pbox3}
s_p(\k)=\sum_{v=0}^{a-1}\dk_{1,v}\geq
\left\lfloor\frac{(m-1)p}{d}\right\rfloor+(a-1)
\left\lceil\frac{p-1}{d}\right\rceil+
\left\lceil\frac{p-j}{d}\right\rceil+A,
\end{eqnarray}
where $A:= \left\lfloor\frac{(m'-1)p}{d}\right\rfloor
-\left\lceil\frac{m'-1}{d}\right\rceil$.
Using $p>d$ one easily observes that $A\geq0$, and the first
part of the lemma follows from (\ref{eq:pbox3}).

Now suppose $p>2d$, $m=1$, the equality holds in (\ref{eq:pbox3}) and
$A=0$. This can only happen if $m'=1$. It follows by induction
that $\dk_{1,0}=\lceil\frac{p-j}{d}\rceil$
and
$\dk_{1,v}=\left\lceil\frac{p-1}{d}\right\rceil$ for
$1\leq v<a-1$. From the equality in (\ref{eq:pbox3}) it follows that
$\dk_{1,a-1}=\left\lceil\frac{p-1}{d}\right\rceil$.
\end{proof}

\begin{lemma}\label{L:simple}
Let $a$ be a natural number and $p>d$. For a polynomial
$f(x)\in W[x]$ of degree $d$ we have
\begin{eqnarray*}
[f(x)^{\sum_{v=0}^{a-1}\left\lceil\frac{p-1}{d}\right\rceil p^v}]_{p^a-1}
&\equiv &
\prod_{v=0}^{a-1}
[(f(x))^{\left\lceil\frac{p-1}{d}\right\rceil}]^{\sigma^a}_{p-1}
\bmod p.
\end{eqnarray*}
\end{lemma}
\begin{proof}
Write
\begin{eqnarray}\label{E:coefficient}
f(x)^{\sum_{v=0}^{a-1}\lceil\frac{p-1}{d}\rceil p^v}
 = \prod_{v=0}^{a-1}f(x)^{\lceil\frac{p-1}{d}\rceil p^v}
\equiv \prod_{v=0}^{a-1}f^{\sigma^v}(x^{p^v})^{\lceil\frac{p-1}{d}\rceil}
           \bmod p.
\end{eqnarray}
Now we write
\begin{eqnarray}\label{E:x}
x^{p^a-1} =\prod_{v=0}^{a-1}x^{p^v (p-1)}.
\end{eqnarray}
Consider contributions of each factor of the product
of~(\ref{E:coefficient}) in the coefficient of~(\ref{E:x}).
Each $v$-th factor of~(\ref{E:coefficient})
contributes to the coefficients of
$x^{p^v m}$ for some $m$, where
$1\leq m\leq d\lceil\frac{p-1}{d}\rceil<2p-1$.
When $v=0$ then it has to contribute to the coefficient of
$x^{p-1}$. Inductively for each $v=1,\ldots,a-1$
the $v$-th factor contributes precisely
the coefficient to $x^{p^v(p-1)}$. It is
easy to see that
\begin{eqnarray*}
[f^{\sigma^v}(x^{p^v})^{\lceil\frac{p-1}{d}\rceil}]_{p^v(p-1)}
&\equiv&
[f(x)^{\lceil\frac{p-1}{d}\rceil}]^{\sigma^v}_{p-1}\bmod p.
\end{eqnarray*}
Thus our assertion follows.
\end{proof}

\begin{lemma}\label{L:pbox3}
Let $p>d$.
Let $a,m,N$ be natural numbers. Let $i,j$ be as in Lemma~\ref{L:1}.
Then
\begin{eqnarray}\label{eq:pbox6}
\ord_p(C_{mp^a-j}(i,N))&\geq& \left\lceil
\frac{
(a-1)\lceil
\frac{p-1}{d}\rceil+\lceil\frac{p-j}{d}\rceil-i}{p-1}\right\rceil.
\end{eqnarray}
Moreover, for $p>2d$ we have
\begin{eqnarray}\label{eq:pbox7}
\ord_p(C_{mp^a-1}(i,N))&=&
\frac{
a\lceil
\frac{p-1}{d}\rceil-i}{p-1}
\end{eqnarray}
if and only if
$$
\left\{
\begin{array}{l}
m=1;\\
a\lceil\frac{p-1}{d}\rceil \equiv i\bmod p-1;\\
\left[ f(x)^{\lceil\frac{p-1}{d}\rceil}\right]_{p-1} \not\equiv 0\bmod p.
\end{array}
\right.
$$
\end{lemma}
\begin{proof}
Let $\k=\ ^t\!(k_1,\ldots,k_d)\in\K_{mp^a-j}$.
Let $\k':=\ ^t\!(k'_1,\ldots,k'_d)$ where
$k'_\ell=\sum_{v=0}^{a-1}\dk_{\ell,v}p^v$,
then  $\k'\in\K_{r'}[a]$.
Let $r':=\sum_{\ell=1}^dk'_\ell$, write
$r'=m'p^a-j$ for some $m'$. From Lemma~\ref{L:pbox1}, it follows that
\begin{eqnarray}\label{E:pbox7}
s_p(\k)
=\sum_{v\geq 0}\dk_{1,v}
\geq\sum_{v=0}^{a-1}\dk_{1,v}=s_p(\k')
\geq (a-1)\left\lceil\frac{p-1}{d}\right\rceil
     +\left\lceil\frac{p-j}{d}\right\rceil.
\end{eqnarray}
Then by (\ref{E:C-E}) and Lemma~\ref{L:4.3}
one easily verifies that (\ref{eq:pbox6}) holds.

Assume (\ref{eq:pbox7}) holds.
Then there is a $\k$ such that the equality in (\ref{E:pbox7})
holds for $j=1$, which implies that $m=1$, $\dk_{1,v}=0$ for $v\geq a$,
$\dk_{1,v}=\lceil\frac{p-1}{d}\rceil$ for $0\leq v\leq a-1$
by Lemma~\ref{L:pbox1}.
Thus $k_1=\sum_{v=0}^{a-1}\lceil\frac{p-1}{d}\rceil p^v$.
So $s_p(k_1)=a\lceil\frac{p-1}{d}\rceil\equiv i\bmod p-1$.
Those $\k\in\K_{mp^a-1}$
which contribute terms in the sum
(\ref{E:C-E}) with minimal valuation necessarily have
$k_1\equiv i\bmod p-1$.
By the identity
$$C_{p^a-1}(i,N)=\sum_{k_1=0}^{\infty}E_{k_1}(i,N)
\cdot
[f(x)^{k_1}]_{p^a-1},
$$
we have by Lemma~\ref{L:4.3}
\begin{eqnarray*}
\ord_p(C_{p^a-1}(i,N))
&\geq &\ord_p(E_{k_1}(i,N))+\ord_p([f(x)^{k_1}]_{p^a-1})\\
&=&\frac{s_p(k_1)-i}{p-1}+\ord_p([f(x)^{k_1}]_{p^a-1})
\end{eqnarray*}
This is equal to
$\frac{a\lceil\frac{p-1}{d}\rceil-i}{p-1}$ if
and only if $[f(x)^{k_1}]_{p^a-1}\not\equiv 0\bmod p$. By
Lemma~\ref{L:simple} this is equivalent to
$
[f(x)^{\lceil\frac{p-1}{d}\rceil}]_{p-1}
\not\equiv 0\bmod p$.

Conversely, the conditions imply that
the contribution of $\k\in\K_{p^a-1}$ with
$k_1=\sum_{v=0}^{a-1}\lceil\frac{p-1}{d}\rceil p^v$
to $\ord_p(C_{p^a-1}(i,N))$ in (\ref{E:C-E})
has valuation $\frac{a\left\lceil\frac{p-1}{d}\right\rceil-i}{p-1}$.
Contribution from other $\k\in\K_{p^a-1}$
has higher valuation by the above arguments.
Thus $\ord_p(C_{p^a-1}(i,N))=
\frac{a\lceil\frac{p-1}{d}\rceil-i}{p-1}$.
This finishes the proof of this lemma.
\end{proof}

\section{Proof of Theorem~\ref{T:box}}\label{section6}

\begin{proof}[Proof of Theorem 1.1]
It suffices to prove the theorem for the case that $\tf(x)$ has
constant coefficient $\ta_0=0$.  On the one hand,
$[\tf(x)^{\left\lceil\frac{p-1}{d}\right\rceil}]_{p-1}$ is independent
of $\ta_0$; on the other hand, the curves $y^p-y=\tf(x)$ and
$y^p-y=\tf(x)+\ta_0$ are isomorphic over ${\overline\F_p}$ for any
$\ta_0$, and hence have the same Newton polygon.  With the assumption
$\ta_0=0$ we can use the results of Section \ref{section3}.

a)
Set $\lambda_0:= \frac{\lceil\frac{p-1}{d}\rceil}{p-1}$.
By the hypothesis on $d$, $p$ and $g$ it is elementary to check that
for $i$ and $j$ in the range of Lemma~\ref{L:1} we have
$(g-2)\left\lceil\frac{p-1}{d}\right\rceil+\left\lceil\frac{p-j}{d}\right\rceil
\geq i$ thus
$\left\lceil (n+g-2)\lambda_0+\frac{\left\lceil\frac{p-j}{d}\right\rceil
-i}{p-1}\right\rceil\geq
\lceil n\lambda_0 \rceil$
for all $n\geq 1$.
By Lemma~\ref{L:pbox3}, we have
\begin{eqnarray*}
\ord_p(C_{mp^{n+g-1}-j}(i,n+g-2))
\geq\left\lceil (n+g-2)\lambda_0
+\frac{\left\lceil\frac{p-j}{d}\right\rceil
-i}{p-1}\right\rceil\geq\lceil n\lambda_0 \rceil.
\end{eqnarray*}
Thus $\NP_1(X/\F_q)\geq \lambda_0$ by Lemma~\ref{keylemma}.

b) Choose a value of $i$ in the range of Lemma~\ref{L:1} for $j=1$
such that the following congruence has a solution for $a$,
\begin{eqnarray*}
a\left\lceil\frac{p-1}{d}\right\rceil &\equiv & i \bmod p-1.
\end{eqnarray*}
For any integer $n>1$ define
\begin{eqnarray*}
\lambda_n &:=& \frac{(n+g-2)\left\lceil\frac{p-1}{d}\right\rceil-i}{(n-1)(p-1)}.
\end{eqnarray*}
Note that $\lambda_n$ is monotonically decreasing as
a function in $n$,
and it converges to $\lambda_0$ as $n$ approaches $\infty$.
Suppose $\NP_1(X/\F_q)>\lambda_0$, then
there exists a positive integer $n_0$ large enough such that
$\NP_1(X/\F_q)>\lambda_{n_0}$. Choose such an $n_0$, and such that
$a=n_0+g-1$ is a solution to the congruence above and
such that
$\frac{(g-1)\left\lceil\frac{p-1}{d}\right\rceil-i}{(p-1)(n_0-1)}\leq1$.
For all $1\leq n <n_0$ we have
\begin{eqnarray*}\label{E:compare}
\lambda_{n_0}\leq\lambda_{n+1}= \frac{(n+g-1)\lceil\frac{p-1}{d}
\rceil-i}{n(p-1)}.
\end{eqnarray*}
Thus, for all $m\geq 1$ and $1\leq n<n_0$ we have by Lemma~\ref{L:pbox3} that
\begin{eqnarray*}
\ord_p(C_{mp^{n+g-1}-1}(i,n+g-2))
&\geq& \left\lceil\frac{(n+g-1)\lceil\frac{p-1}{d}\rceil-i}
{p-1}\right\rceil\\
&\geq& \lceil n\lambda_{n_0} \rceil.
\end{eqnarray*}
On the other hand, since
\begin{eqnarray*}
0<n_0\lambda_{n_0}-\frac{(n_0+g-1)\left\lceil\frac{p-1}{d}\right\rceil-i}
{p-1}=\frac{(g-1)\left\lceil\frac{p-1}{d}\right\rceil-i}{(p-1)(n_0-1)}
\leq1,
\end{eqnarray*}
by our assumption we have
\begin{eqnarray*}
\lceil n_0\lambda_{n_0}\rceil
&=& \frac{(n_0+g-1)\lceil\frac{p-1}{d}\rceil-i}{p-1} + 1.
\end{eqnarray*}
Hence, for all $m\geq 2$ one has by Lemma~\ref{L:pbox3} that
\begin{eqnarray*}
\ord_p(C_{mp^{n_0+g-1}-1}(i,n_0+g-2))
&\geq&\frac{(n_0+g-1)\left\lceil\frac{p-1}{d}\right\rceil-i}{p-1} + 1\\
&  = & \left\lceil n_0\lambda_{n_0}\right\rceil.
\end{eqnarray*}
Hence, the hypotheses of Lemma~\ref{keylemma} are satisfied.
Again by Lemma~\ref{L:pbox3},
\begin{eqnarray*}
\ord_p(C_{p^{n_0+g-1}-1}(i,n_0+g-2)) &\geq&
\left\lceil n_0\lambda_{n_0}\right\rceil - 1,
\end{eqnarray*}
where the equality holds if and only if
$[f(x)^{\lceil\frac{p-1}{d}\rceil}]_{p-1}
\not\equiv 0\bmod p$.
In this case, we have $\NP_1(X/\F_q)<\lambda_{n_0}$ (by
Lemma~\ref{keylemma}), which
contradicts our
assumption that $\NP_1(X/\F_q)>\lambda_{n_0}$.
Therefore, we have $\NP_1(X/\F_q)=\lambda_0$.
\end{proof}

\section{Ax's version of Warning theorem and its application to
slope estimates over finite fields}\label{section7}

This section is independent of the rest of the paper.  Main goal is to
give a lower bound for $\NP_1(X/\F_q)$ without the assumption $p>d$.
Proposition~\ref{P:wan} is due to Daqing Wan.  To begin, we present a
simple lemma due to the fact that we were not able to locate a
suitable reference.

\begin{lemma}\label{L:delight}
Let $X$ be any curve over $\F_q$ where $q=p^\nu$. Let $\lambda$ be a
rational number with $0\leq \lambda\leq 1/2$.
The following two statements are equivalent
\begin{itemize}
\item[(a)] $p^{\lceil \nu n\lambda \rceil}\mid (\#X(\F_{q^n})-1)$
            for all $n\geq 1$;
\item[(b)] $\NP_1(X/\F_q)\geq \lambda.$
\end{itemize}
\end{lemma}
\begin{proof}
The denominator of the $L$ function of $X$ of (genus $g$) is
$P(T)=\prod_{i=1}^{2g}(1-\pi_i T)$
where $\pi_i$'s are eigenvalues of the Frobenius endomorphism of
$X$ relative to $\F_q$. We consider the $\pi_i$ as elements of
$\overline\Q_p$. Extend the valuation $\ord_p(.)$ to $\overline\Q_p$.

Let $\NP_1(X/\F_q)\geq \lambda$. Then
$\ord_p(\pi_i^n)\geq \nu n\lambda$ for all $n\geq 1$
(see \cite[Lemma 4, Chapter IV]{Koblitz}).
But $$\#X(\F_{q^n}) = 1+q^n-\sum_{i=1}^{2g}\pi_i^n,$$ so
$\ord_p(q^n+1-\#X(\F_{q^n}))\geq \nu n\lambda$.
Thus $p^{\lceil \nu n\lambda\rceil}\mid (q^n+1-\# X(\F_{q^n}))$.

Conversely, suppose $p^{\lceil \nu n\lambda\rceil}\mid (q^n+1-\# X(\F_{q^n}))$
for every $n\geq 1$.
One easily derives from (\ref{E:zeta}) that
\begin{eqnarray*}
\exp\left(\sum_{n=1}^{\infty}(q^n+1-\#X(\F_{q^n}))\frac{T^n}{n}\right)
&=&
\frac{1}{\prod_{i=1}^{2g}(1-\pi_i T)}.
\end{eqnarray*}
Taking natural logarithm  and then derivative at both sides, we get
\begin{eqnarray*}\label{E:zeta2}
\sum_{n=1}^{\infty}(q^n+1-\#X(\F_{q^n}))T^{n-1} =
\sum_{i=1}^{2g}\frac{\pi_i}{1-\pi_i T}.
\end{eqnarray*}
Then the left hand side of this power series
converges
$p$-adically for all $T$ with $\ord_p(T)\geq -\nu\lambda$. Comparing to the
right hand side series we have  $\ord_p(\pi_i)\geq \nu\lambda$
for all $i$. Therefore, $\NP_1(X/\F_q)\geq \lambda$.
\end{proof}

\begin{proposition}
\label{P:wan}
Let $X$ be an Artin-Schreier curve over $\overline\F_p$
given by an equation $y^p - y = \tf(x)$ where
$\tf(x)=x^d+\ta_{d-1}x^{d-1}+\ldots+\ta_1x$ and $p\nmid d$.
Then
\begin{eqnarray*}
\NP_1(X/\F_q)
\geq \frac{1}{\displaystyle\max_{\ta_k\neq 0}s_p(k)}
\geq \frac{1}{d}.
\end{eqnarray*}
\end{proposition}
\begin{proof}
Let $\tf(x)$ be a  polynomial over $\F_q$ with $q=p^\nu$
for some $\nu\in\N$.
For any $n\in\N$, write $r=\nu n$.
Let $\{\alpha_k\}_{k=1,\ldots,r}$ be a basis for the degree $r$ extension
$\F_{p^r}/\F_p$. For any $x\in \F_{p^r}$ write
$x = \sum_{t=1}^{r}x_t\alpha_t$ for some $x_t\in \F_p$.
For any $k\in\N$, take its $p$-adic expansion
$k=\sum_{s=1}^{l} k_s p^s$ with $0\leq  k_s\leq p-1$ and
some $l\in\N$. Then
\begin{eqnarray*}
x^k &=& (\sum_{t=1}^{r}x_t\alpha_t)^{\sum_{s=1}^{l} k_s p^s}\\
    &=& (\sum_{t=1}^{r}x_t\alpha_t)^{k_0}(\sum_{t=1}^{r}x_t\alpha_t^p)^{k_1}
        \cdots (\sum_{t=1}^{r}x_t\alpha_t^{p^l})^{k_l}.
\end{eqnarray*} From this, one observes that $\ta_k x^k$ can be
considered as a polynomial in $x_1,\ldots,x_r$
over $\F_{p^r}$ of total degree
$s_p(k)=k_0+k_1+\ldots+k_l.$ Write $\Tr$ for $\Tr_{\F_{p^r}/\F_p}$,
then $\Tr(\tf(x)) = \sum_{k=1}^{d}\Tr (\ta_kx^k)$
is a polynomial in $x_1,\ldots,x_r$
over $\F_p$ of total degree
$D:=\displaystyle\max_{\ta_k\neq 0} s_p(k).$
Then we observe that
$$\#X(\F_{q^n})-1 = p\cdot\#\{x\in\F_{q^n}\mid\Tr(\tf(x)) = 0\}.$$
On the other hand, Ax's theorem~\cite{Ax:64} indicates
$$p^{\lceil\frac{r}{D}\rceil-1}\quad |\quad \#\{x\in\F_{p^r}\mid
\Tr(\tf(x))=0\}.$$
Thus
$$p^{\lceil\frac{\nu n}{D}\rceil}\quad \mid \quad(\#X(\F_{q^n})-1).$$
Applying Lemma~\ref{L:delight}, we have $\NP_1(X/\F_q)\geq \frac{1}{D}$.
The second inequality is elementary.
\end{proof}

\begin{remark}
Note that $\frac{\lceil\frac{p-1}{d}\rceil}{p-1} \geq \frac{1}{d}$
and the equality holds if $p\equiv 1\bmod d$.
Thus for $p>d$
Theorem~\ref{T:box} a) is stronger
than Proposition~\ref{P:wan}.
\end{remark}

\begin{remark}\label{R2:wan}
The supersingularity of curves over $\overline\F_p$ of the form
\begin{eqnarray}\label{E:geer}
y^p-y &=& \sum_{\ell\geq 0}\ta_{p^\ell+1} x^{p^\ell+1},
\end{eqnarray}
as in~\cite{Geer:92,Geer:95} follows from Proposition~\ref{P:wan}.
We conjecture that
if $g=\frac{(p-1)p^h}{2}$ for some $h\geq 1$ then $X$ is
supersingular if and only if $X$ has an equation as in~(\ref{E:geer}).
\end{remark}

\section{A conjecture of Wan}\label{section8}

We first introduce a conjecture of Daqing Wan, then link it to
Artin-Schreier curves. We prove Theorem~\ref{T:2}, which clearly
indicates Theorem~\ref{T:1}.

For every integer $\ell\geq 1$ let
$$S_\ell (f) :=
\sum_{x\in\F_{p^\ell}}\zeta_p^{\Tr_{\F_{p^\ell}/\F_p}(f(x))}.$$
The $L$ function of $f(x)\bmod p$ is defined by
$$L(f\bmod p; T) = \exp\left(\sum_{\ell=1}^{\infty}S_\ell
(f)\frac{T^\ell}{\ell}\right).$$ It is a theorem of
Dwork-Bombieri-Grothendieck that $L(f\bmod p;T)=
1+b'_1T+\ldots+b'_{d-1}T^{d-1}\in\Z[\zeta_p][T]$ for some $p$-th root
of unity $\zeta_p$ in $\overline\Q$. Define the {\em Newton polygon}
of $f\bmod p$ as the lower convex hull of the points
$(k,\ord_p(b'_k))$ in $\R^2$ for $0\leq k\leq d-1$. We denote it by
$\NP(f\bmod p)$. It is the $p$-adic Newton polygon of the polynomial
$L(f\bmod p;T)$.  Define the {\em Hodge polygon} $\HP(f)$ as the
convex hull in $\R^2$ of the points
$(k,\frac{1}{d}+\frac{2}{d}+\ldots+\frac{k}{d})$ for $0\leq k\leq
d-1$.  It is proved by Bombieri~\cite{Bombieri} that the Newton
polygon is always lying above or equal to the Hodge polygon.  See
also~\cite{Sperber} and~\cite{Adolphson-Sperber2} for generalizations.

\begin{remark}
Some literature call
$f\bmod p$ {\em ordinary} if these two polygons coincide~(see~\cite{Katz}).
\end{remark}

\begin{conjecture}[Wan]\label{wan}
There is a Zariski dense subset $\cU$ in $\A^d$ such that for all
$f(x)\in\cU$ we have the following limit exists and
$$\lim_{p\rightarrow \infty}\NP(f\bmod p) = \HP(f).$$
\end{conjecture}

This conjecture was proposed by Wan
in the Berkeley number theory seminar in
the fall of 2000, a general form of which will appear in Section
2.5~\cite{Wan}.  The cases $\deg(f)=3$ and 4 are proved
in~\cite{Sperber} and~\cite{Hong}, respectively. It is also known for
all prime $p\equiv 1\bmod d$, in which case the Newton polygon is
always equal to the Hodge polygon (see~\cite{Adolphson-Sperber}).

It is not hard to verify
\begin{eqnarray*}
L(X/\F_p;T)=\frac{1}{N_{\Q(\zeta_p)/\Q}(L(f\bmod p;T))},
\end{eqnarray*}
where $N_{\Q(\zeta_p)/\Q}(\cdot)$ represents the norm map.
If we ``normalize'' the Newton polygon $\NP(X/\F_p)$  by
shrinking by a factor of $\frac{1}{p-1}$ horizontally and vertically, then
we obtain the Newton polygon $\NP(f\bmod p)$.
Obviously these two Newton polygons have the same ``shape''
thus they have the same slope.
Therefore, Theorem~\ref{T:1} confirms a first slope version of
Conjecture~\ref{wan}.

\begin{theorem}\label{T:2}
Let $d\geq 2$. Let $\cU$ be the set of all monic polynomials $f(x)\in
\A^{d}$ such that $f(x)$ has $$[f(x)^{\left\lceil\frac{p-1}{d}\right
\rceil}]_{p-1}\equiv 0\bmod p$$ for finitely many primes $p$. It is
Zariski dense in $\A^{d}$.  For every $f(x)\in \cU$ we have
\begin{eqnarray*}
\lim_{p\rightarrow \infty} \NP_1(X/\F_p) = \frac{1}{d}.
\end{eqnarray*}
\end{theorem}
\begin{proof}
We fix a natural number $d\geq 2$.  Define
$F(x)=x^d+A_{d-1}x^{d-1}+\ldots+A_0$ as an element of the polynomial
ring $\Q[A_0,\ldots,A_{d-1},x]$ in $d+1$ variables.

Let $k$ be any integer with $0\leq k\leq d-1$ and $\gcd(k-1,d)=1$.
For every prime $p>d$ and $p\equiv 1-k\bmod d$, we write $p=Nd-k+1$
for some integer $N$. So $N=\lceil\frac{p-1}{d}\rceil$. The following
lemmas are suggested by Bjorn Poonen.

\begin{lemma}\label{L:Hpol}
We have $[F(x)^N]_{p-1}\in\Q[A_0,A_1,\ldots,A_{d-1}],$ and it can be
written as a polynomial in $N,A_0,\ldots,A_{d-1}$ with rational
coefficients. Let $f_k(A_0,\ldots,A_{d-1})$ denote the evaluation of
this polynomial at $N=\frac{k-1}{d}$.  Then $f_k(A_0,\ldots,A_{d-1})$
is not the zero polynomial in $\Q[A_0,\ldots,A_{d-1}]$. For any
$f(x)=x^{d}+a_{d-1}x^{d-1}+\ldots +a_0$ in $\A^d$ we have
\begin{eqnarray*}
[f(x)^{\lceil\frac{p-1}{d}\rceil}]_{p-1}
&\equiv&
f_k(a_0,\ldots,a_{d-1})
\bmod p.
\end{eqnarray*}
\end{lemma}
\begin{proof}
Define an auxiliary polynomial $h(T)=A_{d-1}+A_{d-2}T+\ldots+A_0T^{d-1}$.
Note that $[F(x)^N]_{p-1}$
is equal to the  $T^k$-coefficient of
\begin{eqnarray}\label{E:Tk}
(1+Th(T))^N=\sum_{\ell=0}^k\binom{N}{\ell}(Th(T))^\ell.
\end{eqnarray}
For fixed $\ell$ the binomial coefficient $\binom{N}{\ell}$
is a polynomial in $N$ with rational coefficients.
The first assertion follows.
Consider $[F(x)^N]_{p-1}$ as a polynomial in $A_0,\cdots,A_{d-1}$,
its $A_{d-1}^k$-coefficient is equal to
$\binom{N}{k}$ by inspecting (\ref{E:Tk}).
Note $\binom{\frac{k-1}{d}}{k}\neq 0$, so
$f_k$ is not the zero polynomial in $\Q[A_0,\ldots,A_{d-1}]$.
Since $p=Nd-k+1$, we have
$\lceil\frac{p-1}{d}\rceil\equiv\frac{k-1}{d}\bmod p$; hence
\begin{eqnarray*}
[f(x)^{\lceil\frac{p-1}{d}\rceil}]_{p-1}
&\equiv&
f_k(a_0,\ldots,a_{d-1})
\bmod p.
\end{eqnarray*}
This finishes the proof of Lemma~\ref{L:Hpol}.
\end{proof}

\begin{lemma}\label{L:8.5}
Let $f(x)=x^d+a_{d-1}x^{d-1}+\ldots+a_0$ be in $\A^d$.
The following statements are equivalent
\begin{enumerate}
\item[1)] $[f(x)^{\lceil\frac{p-1}{d}\rceil}]_{p-1}
\equiv 0\bmod p$
for infinitely many primes $p$;
\item[2)] there is a $k$ with
$0\leq k\leq d-1$ and $\gcd(k-1,d)=1$, such that
$[f(x)^{\lceil\frac{p-1}{d}\rceil}]_{p-1}
\equiv 0\bmod p$
for infinitely many primes
$p\equiv 1-k\bmod d$;
\item[3)] there is a $k$ with
$0\leq k\leq d-1$ and $\gcd(k-1,d)=1$,
such that
$f_k(a_0,\ldots,a_{d-1})\equiv 0\bmod p$
for infinitely many prime $p\equiv 1-k\bmod d$;
\item[4)] there is a $k$ with
$0\leq k\leq d-1$ and $\gcd(k-1,d)=1$,
such that
$f_k(a_0,\ldots,a_{d-1})=0$.
\end{enumerate}
\end{lemma}
\begin{proof}
Parts 1) and 2) are clearly equivalent.
Parts 2) and 3) are equivalent by Lemma~\ref{L:Hpol}.
Parts 3) and 4) are equivalent
because $f_k(a_0,\ldots,a_{d-1})\in\Q$
has to vanish if it vanishes modulo $p$ for infinitely many primes $p$;
conversely, since $\gcd(k-1,d)=1$, there are infinitely many prime
$p\equiv 1-k\bmod d$ by Dirichlet. This concludes Lemma~\ref{L:8.5}.
\end{proof}

The complement $\cU^c$ of $\cU$ in $\A^d$ is the set of all $f(x)$ in
$\A^d$ such that $$[f(x)^{\lceil \frac{p-1}{d}\rceil}]_{p-1}\equiv
0\bmod p$$ for infinitely many prime $p$.  Write
$G(A_0,\ldots,A_{d-1}):=\prod_{k}f_k(A_0,\ldots,A_{d-1})$.
By
Lemma~\ref{L:8.5}, $\cU^c$ is equal to the set of all
$f(x)=x^d+a_{d-1}x^{d-1}+\ldots+a_0$ in $\A^d$ such that
$(a_0,\ldots,a_{d-1})$ is a zero of a polynomial $G(A_0,\ldots,A_{d-1})$
in $\Q[A_0,\ldots,A_{d-1}]$.  But
$G(A_0,\ldots,A_{d-1})$ is not the zero polynomial by
Lemma~\ref{L:Hpol}, so $\cU^c$ is Zariski closed. Therefore, $\cU$ is
Zariski open and hence dense in $\A^d$. This concludes the first part
of the theorem.

Now let $f(x)\in\cU$. Then there exists a natural number $M$ such that
for all $p\geq M$ we have
\begin{eqnarray*}
[f(x)^{\left\lceil\frac{p-1}{d}\right\rceil}]_{p-1}
     &\not\equiv& 0\bmod p,
\end{eqnarray*}
and hence $\NP_1(X/\F_p) = \frac{\lceil\frac{p-1}{d}\rceil}{p-1}$
by Theorem~\ref{T:box}.
Therefore, for every $f(x)\in\cU$ we have
\begin{eqnarray*}
\lim_{p\rightarrow \infty}\NP_1(X/\F_p) &=& \frac{1}{d}.
\end{eqnarray*}
This finishes the proof of Theorem~\ref{T:2}.
\end{proof}

\end{document}